\documentclass[leqno,12pt]{amsart}\usepackage{amsmath,amstext,amssymb,amsopn,amsthm,mathrsfs}

\theoremstyle{plain}
\newtheorem{theorem}{Theorem}[section]
\newtheorem{defi}{Definition}[section]

\newtheorem{lemma}{Lemma}[section]

\numberwithin{equation}{section}
\newcommand{\dem}{\medskip \par \noindent \mbox{\bf Proof. }}
\def\ep{\hfill{$\Box $}}

\newcommand{\R}{{\mathcal R}}
\newcommand{\be}{\begin{equation}}
\newcommand{\ee}{\end{equation}}

\begin{document}

\title{On the Theta semigroup}
\author{Wilfredo~O.~Urbina }
\address{Department of Mathematics and Actuarial Sciences, Roosevelt University, Chicago, IL,
   60605, USA.}
\email{wurbinaromero@roosevelt.edu}

\author{Ahmed Zayed}
\address{Department of Mathematical Sciences, DePaul University, Chicago, IL,
   60614, USA.}
   \email{azayed@depaul.edu}



\begin{abstract}
In this paper we consider a semigroup on trigonometric expansions  $\{  {\mathcal T}_t: t \geq 0\}$  that will be called the Theta semigroup since  its kernel is a multiple of $\theta_3(x,q)$, the third Jacobi theta function. We study  properties of this semigroup and  prove that  it is a positive diffusion semigroup. We also obtain that its subordinated semigroup  is the classical Poisson semigroup. The extensions to higher dimensions and to periodic ultra distributions  are also considered.
\end{abstract}

\maketitle

\section{Introduction}
There are two classical semigroups in analysis in $\mathbb{R}^n$, the heat semigroup,
\begin{eqnarray*}\label{semicalor}
T_t f(x)  = \frac{1}{(4\pi t)^{n/2}}\int_{\mathbb R^n} e^{\frac{-|x-y|^2}{4t}} f(y)\,dy
 \end{eqnarray*}
 and the Poisson semigroup,
 \begin{eqnarray*}
P_t f(x) &=&  \frac{\Gamma(\frac{n+1}{2})}{ \pi^{(n+1)/2}} \int_{\mathbb R^n} \frac{t}{ (|x-y|^2+t^2)^{(n+1)/2}} f(y)\,dy.
\end{eqnarray*}
They are related by Bochner's subordination formula, see E. Stein \cite{st1}
 \begin{equation}\label{bochner}
 e^{-\lambda} = \frac{1}{\sqrt \pi} \int_0^{\infty} \frac{e^{-u}}{\sqrt u} e^{-\lambda^2/4u} du,
\end{equation}
since
\begin{eqnarray*}
P_t f(x) &=& \frac{1}{\sqrt \pi} \int_0^{\infty} \frac{e^{-u}}{\sqrt u} T_{t^2/4u}f(x) du.
\end{eqnarray*}

In the periodic case the situation is a little different. Let ${\mathbb T}$ be the unit circle and consider the family of trigonometric functions ${\mathcal E}=\{e^{nix}\}$, which is a complete orthogonal system in $L^2({\mathbb T})$. Given $f \in L^2({\mathbb T})$, let us consider its Fourier coefficients
\begin{equation}\label{coef}
\hat{f}(n) =  \frac{1}{2\pi}\int_{\mathbb T} f(y) e^{-niy} dy,
\end{equation}
then we have that the Fourier series associated to $f$ is given by
\begin{equation}\label{series}
 \sum_{n=-\infty}^{\infty}   \hat{f}(n)  e^{nix}.
\end{equation}
It is well known that the partial sums of the  Fourier series tend to $f$ in $L^2$-norm.
Moreover for $f \in L^1({\mathbb T})$, we can associate a Fourier series to $f$ using the same argument, but in that case the convergence of the series is
a more difficult problem.\\

In $\mathbb T$ the Poisson integral is defined by means of the Abel's summability of the Fourier series of $f$,
\begin{eqnarray*}\label{seriesAbel}
 \sum_{n=-\infty}^{\infty}  r^{|n|} \hat{f}(n)  e^{nix} &=& \frac{1}{2\pi}  \int_{\mathbb T} [\sum_{n=-\infty}^{\infty}  r^{|n|}  e^{ni(x-y)} ] f(y) dy \\
 &=& \frac{1}{2\pi}  \int_{\mathbb T} \frac{1-r^2}{1-2r \cos(x-y) +r^2} f(y) dy,
\end{eqnarray*}
 since as it is well known
\begin{eqnarray*}
 \sum_{n=-\infty}^{\infty}  r^{|n|}  e^{nix}  =  \frac{1-r^2}{1-2r \cos x +r^2}
\end{eqnarray*}
The classical Poisson semigroup on ${\mathbb T}$ is obtained  by taking $r = e^{-t},$ i.e.
\begin{equation}
{\mathcal P}_t f(x) =   \frac{1}{2\pi}  \int_{\mathbb T}   \frac{1- e^{-2t}}{1- 2 e^{-t} \cos (x-y) +e^{-2t}} f(y) dy.
\end{equation}
$\{{\mathcal P}_t\}$  is an $L^p$-contraction semigroup, $1\leq p \leq \infty$, which is positive and conservative. In the next section we are going to construct another semigroup, that we will call the Theta semigroup $\{{\mathcal T}_t\}$ such that the Poisson semigroup will be its subordinated semigroup by Bochner's subordination formula (\ref{bochner}).\\

 The article is organized as follows. In section 2 we define and study the main properties of the Theta semigroup as well as its generalization to $d$ dimensions. In section 3 we extend the Theta semigroup to ultra-distributions.

\section{Definition and main properties of the Theta semigroup.}

For each $t \geq 0$ and $f \in L^2({\mathbb T})$ let us consider the multiplier operator
$$ {\mathcal T}_t f(x) =  \sum_{n=-\infty}^{\infty}  \hat{f}(n) e^{-n^2t} e^{nix}.$$
Then we have the following integral representation of $ {\mathcal T}_t$,

\begin{eqnarray}
  {\mathcal T}_t f(x) &=&  \sum_{n=-\infty}^{\infty}  \hat{f}(n) e^{-n^2t} e^{nix} = \frac{1}{\pi} \sum_{n=-\infty}^{\infty}  \int_{\mathbb T} f(y) e^{-niy} dy e^{-n^2t} e^{nix} \nonumber \\
 &=& \frac{1}{2\pi}  \int_{\mathbb T}  f(y) [ \sum_{n=-\infty}^{\infty}e^{ni(x-y)} e^{-n^2t} ] dy= \frac{1}{2\pi}  \int_{\mathbb T}  f(y) \theta_3(x-y,e^{-t})dy \nonumber\\
&=& \int_{\mathbb T}  f(y)K_t(x-y)dy = (f \ast K_t)(x), \label{Thetatransf}
\end{eqnarray}
where
$$ K_t(x) =  \frac{1}{2\pi}  \theta_3(x, e^{-t}),$$
and
\begin{eqnarray*}
\theta_3(x,q) &=& \sum_{n=-\infty}^{\infty} q^{n^2} e^{nix} = 1+ 2\sum_{n=1}^{\infty} q^{n^2} \cos (nx)\\
&=& \prod_{n=1}^{\infty} [ 1 + 2 q^{2n-1} \cos x + q^{2(2n-1)}](1-q^{2n}),
\end{eqnarray*}
is the third Jacobi theta function, $|q| <1$. We will call the the integral transform in (\ref{Thetatransf}) the theta transform.

Using this integral representation and the trivial estimate

\begin{equation}\label{thetabound}
|\theta_3(x,q)| \leq \ 1+ 2\sum_{n=1}^{\infty} q^{n^2} = \theta_3(0,q),
\end{equation}
 we can extend $\{\mathcal{T}_t\}$ to $L^1([0,2\pi])$, since
$$
|| {\mathcal T}_t f ||_1\leq \int_{\mathbb T}  \int_{\mathbb T} |f(y)| |\theta_3(x-y,e^{-t})| dy dx
\leq  \theta_3(0,e^{-t}) \int_{\mathbb T} |f(y)| dy = C  ||f ||_1,
$$
 with $C=  \theta_3(0,e^{-t}) \neq 0$. \\

For $\{  {\mathcal T}_t: t \geq 0\}$ we have,

 \begin{theorem}
 The family of operators $\{  {\mathcal T}_t: t \geq 0\}$ satisfies the following properties:
\begin{enumerate}
\item[i)] Semigroup property:
 For any $t_1, t_2 \geq 0,  {\mathcal T}_{t_1 + t_2} =  {\mathcal T}_{t_1} \circ  {\mathcal T}_{t_2}.$
\item[ii)] Positivity and conservative property:
$ {\mathcal T}_t 1 = 1$ and if $ f\geq 0$ then  $ {\mathcal T}_t f \geq 0$ for all $t \geq 0$.
\item[iii)] Contractivity property:
 For any $ t \geq 0$ and $1 \leq  p \leq \infty $,
$$||  {\mathcal T}_tf ||_{p} \leq ||f||_{p}.$$
\item[iv)]  Strong $L^p$-continuity property:
For any $1 \leq p < \infty $ and all $ f \in  L^p({\mathbb T})$
the mapping $ t \rightarrow  {\mathcal T}_tf$ is continuous from $[0,\infty)$ to $L^p({\mathbb T}).$
\item[v)] Symmetry property: for any $ t\geq 0,   {\mathcal T}_t$ is a self-adjoint operator in  $ L^2({\mathbb T})$:
\begin{equation}
\int_{\mathbb T} {\mathcal T}_tf(x) g(x) dx = \int_{\mathbb T}f(x)
 {\mathcal T}_tg(x) dx.
\end{equation}
In particular, the Lebesgue measure is the invariant measure for $\{ {\mathcal T}_t\}$,
\begin{equation}
\int_{\mathbb T}  {\mathcal T}_tf(x) dx = \int_{\mathbb T}f(x) dx,
\end{equation}
\item[vi)] Infinitesimal generator: $L = \frac{d^2}{dx^2} $  is the infinitesimal generator
 of $\{  {\mathcal T}_t: t \geq 0\}$,
\begin{equation}
 \lim_{t \rightarrow 0} \frac{ {\mathcal T}_tf - f}{t} = Lf.
 \end{equation}
 \item[vii)] $u(x,t)  =  {\mathcal T}_t f(x)$ is solution of the heat equation on ${\mathbb T} \times\mathbb{R}^+$,
  $$ \frac {\partial  u(x,t) }{\partial t}  =\frac {\partial^2 u(x,t) }{\partial x^2}, \quad (x,t) \in {\mathbb T} \times\mathbb{R}^+,$$
  with initial condition $u(x,0)  = f(x), \, x \in{\mathbb T},$  and non-homogeneous boundary conditions.
\end{enumerate}
\end{theorem}
\dem
\begin{enumerate}
\item [i)] By definition it is clear that $\{ {\mathcal T}_t\}$ is a semigroup in $L^2({\mathbb T})$, since
\begin{equation*}
 {\mathcal T}_{t_1}({\mathcal T}_{t_2} f) =
    {\mathcal T}_{t_1}(\sum_{n=-\infty}^{\infty}  \hat{f}(n) e^{-n^2t_2} e^{nix} )=  \sum_{n=-\infty}^{\infty}  \hat{f}(n)  e^{-n^2(t_1+t_2)} e^{nix}= {\mathcal T}_{t_1+t_2} f.
\end{equation*}
For the case $p \neq 2$,  we have that, by the orthogonality of ${\mathcal E}=\{e^{nix}\}$, the kernel $K_t$ satisfies the Chapman-Kolmogorov equation,
\begin{equation}
\int_{\mathbb T} K_{t_2} (y) K_{t_1}(x-y) dy  = K_{t_1 + t_2} (x).
\end{equation}
Therefore, by Fubini's theorem
\begin{eqnarray*}
 {\mathcal T}_{t_1}({\mathcal T}_{t_2} f(x)) 
  &=& \int_{\mathbb T} f(u) [\int_{\mathbb T} K_{t_2}(y-u) K_{t_1}(x-y)dy]  du \\
   &=& \int_{\mathbb T} f(u)  K_{t_1+t_2}(x-u) du = {\mathcal T}_{t_1+t_2} f(x).
\end{eqnarray*}

\item[ii)] Let us observe that
\begin{eqnarray}\label{consevative}
\nonumber \int_{\mathbb T}  K_t(x) dx &=&  \frac{1}{2\pi}  \int_{\mathbb T}   \theta_3(x,e^{-t}) dx =   \frac{1}{2\pi}  \int_{\mathbb T} (1+ 2 \sum_{n=1}^{\infty} e^{-t n^2} \cos (nx)) dx\\
 &=&1+  \frac{1}{\pi}  \sum_{n=1}^{\infty} e^{-t n^2}  \int_{\mathbb T} \cos (nx) dx =1,
\end{eqnarray}
and therefore  trivially we get that $ {\mathcal T}_t 1 = 1$ i.e. $\{ {\mathcal T}_t\}$ is a conservative semigroup. Now we have to prove that $ {\mathcal T}_t $ preserves positivity. This will follows immediately from the fact that $\theta_3(x,q)$ is a positive function, which can be obtained immediately from the product representation of   $\theta_3(x,q)$:
 \begin{eqnarray*}
\theta_3(x,q) &=& \prod_{n=1}^{\infty} [ 1 + 2 q^{2n-1} \cos x + q^{2(2n-1)}](1-q^{2n})\\
\end{eqnarray*}
 observing that each factor $(1 + 2 q^{2n-1} \cos x + q^{2(2n-1)})$ and  $(1-q^{2n})$ are non
negative for real $q$ with $|q| <1$. So $\theta_3(x,q) \geq 0$.
 
\item[iii)] The proof  that  ${\mathcal T}_t$ is a $L^p({\mathbb T})$-contraction, i.e.,
\begin{equation}\label{lpcont}
 || {\mathcal T}_t f ||_p \leq  ||f ||_p
\end{equation}
for any $f \in L^p({\mathbb T})$, $1 < p < \infty,$ is trivial since by Young's inequality
 \begin{eqnarray*}
|| {\mathcal T}_t f ||_p =  ||f \ast K_t ||_p\leq  ||f ||_p || K_t ||_1 = ||f ||_p.
\end{eqnarray*}

\item[iv)] The fact that the Theta semigroup is $L^p$-strongly continuous for $1 \leq p \leq \infty$, means that the map
$ {\mathcal T}_t \rightarrow  {\mathcal T}_{t_0}$ in $L^p$-norm as $t \rightarrow t_0$. By the semigroup property it is enough to prove $ {\mathcal T}_t \rightarrow  {\mathcal T}_{0}$ in $L^p$-norm as $t \rightarrow 0$.
First, observe that, by Parseval's identity, the case $p=2$ is trivial, since
\begin{eqnarray*}
\| {\mathcal T}_t f  - f \|_2^2&=& \sum_{n=-\infty}^{\infty} | \hat{f}(n)|^2 |e^{-n^2t} -1|^2
\end{eqnarray*}
and the last term tend to zero as $t \rightarrow 0$. Now if $p \neq 2$, we have
\begin{eqnarray*}
| {\mathcal T}_tf(x) - f(x)| &\leq& \int_{\mathbb T}  |f(y) - f(x)| K_t(x-y) dy \\
&=& \int_{|x-y| < \delta}  + \int_{|x-y| \geq \delta}  |f(y) - f(x)| K_t(x-y) dy\\
&=& (I) + (II)
\end{eqnarray*}
It is enough to prove it for a dense class. Let $f \in \mathcal{P}$ be a polynomial, then given $\epsilon >0$, let us take $\delta > 0$ such that
$ |f(y) - f(x)| < \epsilon$ if $|x-y| < \delta$. Therefore
we get that
$$ |(I)| <\epsilon  \int_{|x-y| < \delta} K_t(x-y) dy \leq \epsilon.$$
For the second integral we have trivially
$$|(II)| < 2 \| f\|_{\infty} \int_{|x-y| \geq \delta} K_t(x-y) dy,$$
and $ \int_{|x-y| \geq \delta} K_t(x-y) dy \rightarrow 0$ as $t \rightarrow 0$, since the theta function tends to the Dirac comb $\sum_{n=-\infty}^{\infty} \delta(x-n)$ as $t \rightarrow 0$.

\item[v)]  The invariant measure for $\{ {\mathcal T}_t\}$ is the Lebesgue measure  on ${\mathbb T}$  since
 $$  \int_{\mathbb T}   {\mathcal T}_tf(y) dy =  \int_{\mathbb T} [  \sum_{n=-\infty}^{\infty}  \hat{f}(n) e^{-n^2t} e^{niy} ] dy =2\pi \hat{f}(0)= \int_{\mathbb T} f(y)  dy.$$
Moreover, since the theta function is an even function of $x,$ then it is easy to see that $\{ {\mathcal T}_t\}$ is a symmetric semigroup
$$  \int_{\mathbb T}   {\mathcal T}_tf(x) g(x) dx =  \int_{\mathbb T}  f(y)  {\mathcal T}_tg(y) dy.$$
Therefore $T_t$ is (formally) self-adjoint.

\item[vi)] The infinitesimal generator  $L$ of $\{ {\mathcal T}_t\}$ is by definition
\begin{eqnarray*}
 {\mathcal T} f(x)&=& \lim_{t \rightarrow 0} \frac{ {\mathcal T}_t f(x) -f(x)}{t} =  \sum_{n=-\infty}^{\infty}  \hat{f}(n)\lim_{t \rightarrow 0} \frac{ e^{-n^2t} -1}{t} e^{nix}.\\
& =&-  \sum_{n=-\infty}^{\infty}  n^2 \hat{f}(n)  e^{-n^2t} e^{nix} =  \frac{d^2}{dx^2} (  \sum_{n=-\infty}^{\infty} \hat{f}(n)  e^{-n^2t} e^{nix})\\
&=& \frac{d^2 f(x)}{dx^2} .
\end{eqnarray*}
 Thus the infinitesimal generator is $L = \ \frac{d^2}{dx^2} $ and
 $$ D(L) = \{ f\in L^2({\mathbb T}): \sum_{n=-\infty}^{\infty}  n^4 |\hat{f}(n)|^2 < \infty\}.$$
\item[vi)]  If $q= e^{\pi \tau i}, \Re(\tau) >0 $ then the third Jacobi theta function  $\theta_3(x, \tau)$ satisfies the differential equation
 $$  i \pi \frac {\partial^2  \theta_3(x, \tau)}{\partial x^2} +  \frac {\partial  \theta_3(x, \tau)}{\partial \tau} =0,$$
 since
 $$ \frac {\partial^2  \theta_3(x, \tau)}{\partial x^2} =-  \sum_{n=-\infty}^{\infty} n^2 e^{\pi \tau i n^2} e^{nix}  \, \; \mbox{and} \; \frac {\partial  \theta_3(x, \tau)}{\partial \tau} = i \pi  \sum_{n=-\infty}^{\infty} n^2 e^{\pi \tau i n^2} e^{nix}.$$
 In our case $ q = e^{\pi (i t/\pi) i}$ i.e. $\tau =i t/\pi$ which clearly satisfies $\Re(\tau) >0 $, and the differential equation becomes
  $$ \frac {\partial^2  \theta_3(x, \tau)}{\partial x^2} =  \frac {\partial  \theta_3(x, t)}{\partial t}.$$
  Therefore $\{ {\mathcal T}_t \}$ satisfies
  $$ L {\mathcal T}_t f(x) =  \frac {\partial^2 T_t f(x)}{\partial x^2} = \frac {\partial  T_t f(x)}{\partial t},$$
  i.e., $u(x,t)  =  {\mathcal T}_t f(x)$ is a solution of the heat equation on $[0,2\pi] \times\mathbb{R}^+$,
  $$  L u(x,t)  = \frac {\partial^2 u(x,t) }{\partial x^2} = \frac {\partial  u(x,t) }{\partial t}, \quad (x,t) \in [0,2\pi] \times\mathbb{R}^+,$$
  with initial condition $u(x,0)  = f(x), \, x \in[0,2\pi],$ and non-homogeneous boundary conditions
  at $x=0$ and $x=2\pi$ ($u(x,0) = u(2\pi, t) = \sum_{n=0}^{\infty} \hat{f}(n) e^{-n^2t}).$ Thus  $u(x,t)$ is the parabolic extension of $f$ to the
  strip $[0,2\pi] \times\mathbb{R}^+.$ \ep

\end{enumerate}

We will call  $\{ {\mathcal T}_t\}$ the {\bf Theta semigroup}, it is a positive symmetric diffusion semigroup, see \cite{st2} and it can be extended analytically to the sector $S_p: | \arg(t+i \tau) < \pi/2(1-|2/p -1|)$, $1 < p < \infty$.

Moreover,  we shall see in the next section that the Theta semigroup $\{ {\mathcal T}_t\}$ can be extended to a class of periodic ultra distributions.\\

 Considering the maximal function of the semigroup
\begin{equation}
 {\mathcal T}^* f(x) = \sup_{t>0} | {\mathcal T}_t f(x)|,
\end{equation}
for any $f \in L^p([0, 2\pi]),$ and using the maximal theorem \cite[p. 73]{st2}, we have
\begin{equation}
\|  {\mathcal T}^* f\|_p \leq \| f\|_p, \quad 1 < p  \leq \infty,
\end{equation}
and
\begin{equation}
\lim_{t \rightarrow 0} {\mathcal T}_t f(x) = f(x), \, \mbox{a.e.} \quad 1 < p  < \infty.
\end{equation}

 Now using the {\em Bochner's subordination formula},
\begin{equation}\label{bochnerf}
 e^{-\lambda} = \frac{1}{\sqrt \pi} \int_0^{\infty} \frac{e^{-u}}{\sqrt u} e^{-\lambda^2/4u} du, \,\lambda >0,
\end{equation}
 we define the subordinated semigroup as
\begin{equation}
{\mathcal  P}_t f(x) = \frac{1}{\sqrt \pi} \int_0^{\infty} \frac{e^{-u}}{\sqrt u}  {\mathcal T}_{(t^2/4u)}f(x) du.
\end{equation}
Then
\begin{eqnarray*}
{\mathcal P}_t f(x) &=&  \frac{1}{\sqrt \pi} \int_0^{\infty} \frac{e^{-u}}{\sqrt u} \frac{1}{2\pi}  \int_{\mathbb T}  f(y) [ \sum_{n=-\infty}^{\infty}e^{ni(x-y)} e^{-n^2t^2/4u} ] dy du\\
&=&  \frac{1}{2\pi}  \int_{\mathbb T}  f(y) [ \sum_{n=-\infty}^{\infty}e^{ni(x-y)}  \frac{1}{\sqrt \pi} \int_0^{\infty} \frac{e^{-u}}{\sqrt u} e^{-n^2t^2/4u}du  ] dy \\
&=&  \frac{1}{2\pi}  \int_{\mathbb T}  f(y) [ \sum_{n=-\infty}^{\infty}e^{ni(x-y)}  e^{-|n|t}  ] dy,
\end{eqnarray*}
and since
\begin{eqnarray*}
 \sum_{n=-\infty}^{\infty}e^{nix}  e^{-|n|t} =\frac{1- e^{-2t}}{1- 2 e^{-t} \cos x+e^{-2t}},\
\end{eqnarray*}
we have
$$ {\mathcal P}_t f(x) =   \frac{1}{2\pi}  \int_{\mathbb T} \frac{1- e^{-2t}}{1- 2 e^{-t} \cos (x-y) +e^{-2t}}  f(y)dy.$$

Therefore,  the subordinate semigroup of the Theta semigroup $\{ {\mathcal T}_t\}$ is the classical Poisson semigroup $\{{\mathcal P}_t\}$ on ${\mathbb T}$.\\

Now, by the tensorization argument we can extend the Theta semigroup  to  higher dimensions. That is to say, considering $f \in L^1({\mathbb T}^d)$, we can define the $d$-dimensional Theta semigroup as
\begin{equation}
 {\mathcal T}^d_t f(x) = \int_{{\mathbb T}^d}  f(y)K_t(x-y)dy = (f \ast K_t)(x),
\end{equation}
where
\begin{equation}
K_t(x) =  \frac{1}{(2\pi)^d}  \prod_{i=1}^d\theta_3(x_i, e^{-t}),\quad x = (x_1, \cdots, x_d).
\end{equation}
Let us observe that by (\ref{thetabound}) we have
$$ |K_t(x)| \leq  \prod_{i=1}^d\theta_3(0, e^{-t}) =  \theta_3(0, e^{-t})^d,$$
and, therefore,  $ {\mathcal T}^d_t f$ is well defined for $f \in L^1({\mathbb T}^d).$ Moreover, we have
\begin{theorem}\label{thetaprop}
 The family of operators $\{  {\mathcal T}^d_t: t \geq 0\}$ satisfies the following properties:
\begin{enumerate}
\item[i)] Semigroup property:
 For any $t_1, t_2 \geq 0,  {\mathcal T}^d_{t_1 + t_2} = T^d_{t_1} \circ T^d_{t_2}.$
\item[ii)] Positivity and conservative property:
$ {\mathcal T}^d_t 1 = 1$ and if $ f\geq 0$ then  $ {\mathcal T}^d_tf \geq 0$ for all $t \geq 0$.
\item[iii)] Contractivity property:
 For any $ t \geq 0$ and $1 \leq  p \leq \infty $,
$$||  {\mathcal T}^d_tf ||_{p} \leq ||f||_{p}.$$
\item[iv)]  Strong $L^p$-continuity property:
For any $1 \leq p < \infty $ and all $ f \in  L^p({\mathbb T}^d)$
the application $ t \rightarrow  {\mathcal T}^d_tf$ is continuous from $[0,\infty)$ to $L^p({\mathbb T}^d).$
\item[v)] Symmetry property: for any $ t\geq 0,   {\mathcal T}^d_t$ is a self-adjoint operator in  $ L^2({\mathbb T}^d)$:
\begin{equation} \label{autoadOU}
\int_{{\mathbb T}^d}  {\mathcal T}^d_tf(x) g(x) dx = \int_{{\mathbb T}^d}f(x)
 {\mathcal T}^d_tg(x) dx.
\end{equation}
In particular, the Lebesgue measure is the invariant measure for $\{ {\mathcal T}^d_t\}$,
\begin{equation}
\int_{{\mathbb T}^d}  {\mathcal T}^d_tf(x) dx = \int_{{\mathbb T}^d}f(x) dx,
\end{equation}
\item[vi)] Infinitesimal generator: $L =  \sum_{i=1}^d \frac{\partial^2 }{\partial x_i^2}=\Delta $  is the infinitesimal generator
of $\{  {\mathcal T}^d_t: t \geq 0\}$,
\begin{equation}
 \lim_{t \rightarrow 0} \frac{{\mathcal T}^d_tf - f}{t} = Lf.
\end{equation}
 \item[vii)] $u(x,t)  =  {\mathcal T}_t f(x)$ is solution of the heat equation on $[0,2\pi]^d \times\mathbb{R}^+$,
  $$  \frac {\partial  u(x,t) }{\partial t}  = \Delta u(x,t) , \quad (x,t) \in [0,2\pi]^d \times\mathbb{R}^+,$$
\end{enumerate}
\end{theorem}
\dem

Observe that
\begin{eqnarray*}
 {\mathcal T}^d_t f(x)  &=& \frac{1}{(2\pi)^d}  \int_{{\mathbb T}^d}  f(y) \prod_{i=1}^d\theta_3(x_i-y_i, e^{-t})dy_1\cdots dy_d\\
  &=& \frac{1}{(2\pi)^d}  \int_{\mathbb T} \cdots (\int_{\mathbb T}  f(y_1,\cdots, y_d) \theta_3(x_1-y_1, e^{-t})dy_1 )\cdots   \theta_3(x_d-y_d, e^{-t}) dy_d,\\
\end{eqnarray*}
thus the results follows from Theorem 1 by iterating the one variable argument $d$-times and using Fubini's theorem and Minkowski integral inequality for property iii).\ep

Finally, if we consider the subordinated semigroup $\{{\mathcal P}^d_t\}$ of the Theta semigroup$\{{\mathcal T}^d_t\}$, using Bochner subordination formula (\ref{bochnerf}), we get,
\begin{eqnarray*}
{\mathcal P}^d_t f(x) &=& \frac{1}{\sqrt \pi} \int_0^{\infty} \frac{e^{-u}}{\sqrt u}  {\mathcal T}^d_{(t^2/4u)}f(x) du\\
&=&  \frac{1}{\sqrt \pi} \int_0^{\infty} \frac{e^{-u}}{\sqrt u} \frac{1}{(2\pi)^d}  \int_{{\mathbb T}^d} f(y) [\prod_{j=1}^d \sum_{n_j=-\infty}^{\infty}e^{n_ji(x_j-y_j)}  e^{-n_j^2t^2/4u} ] dy du\\
&=&  \frac{1}{(2\pi)^d}  \int_{{\mathbb T}^d}f(y)[ \sum_{n_1,\cdots , n_d=-\infty}^{\infty}e^{\sum_{j=1}^d n_ji(x_j-y_j)}\\
&\times & \quad \quad \quad  \frac{1}{\sqrt \pi} \int_0^{\infty} \frac{e^{-u}}{\sqrt u} e^{- (\sum_{j=1}^d n_j^2) t^2/4u}du]dy \\
&=&  \frac{1}{(2\pi)^d}  \int_{{\mathbb T}^d}  f(y) [ \sum_{n_1,\cdots , n_d=-\infty}^{\infty}e^{\sum_{j=1}^d n_ji(x_j-y_j)}  e^{-t \sqrt{\sum_{j=1}^d n_j^2} }  ] dy. \\
\end{eqnarray*}
Let us observe that $\{{\mathcal P}^d_t\}$ is the (classical) Poisson semigroup only in the case that $d=1$.

\section{Theta Semigroup on a Space of Ultra-distributions}

 In this section we extend the theta transform as well as the structure of the theta semigroup to a space of ultra-distributions
 on the unit circle. This space of ultra-distributions contains the space of hyperfunctions on the unit circle as a
 proper subspace. A hyperfunction on the unit circle can be viewed as an element of the dual space of the space of
 analytic functions on the unit circle or as the weak limit of a harmonic function $h(r, \theta)$ in the unit disc as $r$
 approaches one ($r\rightarrow 1^-$). We employ the convolution structure on the unit circle together with the theta transform to obtain a
 weak solution of the heat equation on ${\mathbb T}\times\mathbb{R}^+$ in the sense of ultra-distributions.

 The results of this section are parallel to those obtained by Betancor et al in a series of papers \cite{BetanBelhad1,BetanBelhad2,Betan3}
 to extend the
 Hankel transform and convolution to a space of Beurling distributions and hyperfunctions. The Hankel transform appears naturally
 in solving the initial-value problem related to the generalized heat equation $$\frac{\partial u(x,t)}{\partial t}= L_x u(x,t),\quad u(x,0)=f(x), \quad (x,t)\in \mathbb{R}^+\times\mathbb{R}^+ ,$$
 where $L_x$ is the Bessel differential operator; see \cite{ZayedHaimo1,ZayedHaimo2} for more  general types of generalized heat equations. The space of hyperfunctions considered in \cite{Betan3}
  is the dual space of the space of even entire functions $f(x)$ satisfying
  $$\sup_{x\in\R^+,k\in\mathbb{N}}\frac{\left|f^{(k)}(x)\right|e^{m|x|}}{h^k k!},$$
  for every $m\in \mathbb{N}$ and $h>0.$

Let us consider $C^\infty ({\mathbb T})$ the set of all infinitely
differentiable functions
 on ${\mathbb T}$ or equivalently the set of all  infinitely differentiable functions that
 are periodic with period $2\pi.$
 We endow $C^\infty ({\mathbb T})$ with the topology of uniform convergence so that
  $\left\{f_n\right\}_{n=0}^\infty \subset C^\infty ({\mathbb T})$ converges to $f$
   in this topology if and only if $f^{(k)}_n\rightarrow f^{(k)}$
    uniformly on ${\mathbb T}$, for all $k=0,1, \cdots ,$ as $n\rightarrow \infty.$

 It is known \cite{Schwartz, Zeman} that if $f(t)= \sum_{n=-\infty}^{\infty} \hat{f}_n e^{i n t }\in C^\infty ({\mathbb T}),$
 then the Fourier coefficients  $\hat{f}_n$ of $f$ satisfy the relation
 $$\left( 1+|n|\right)^p |\hat{f}_n| \leq B , \mbox{ for all } 0\leq p \mbox{ and some } 0<B.$$
 The dual space of  $C^\infty ({\mathbb T})$, which we will denote by ${\mathcal E} ({\mathbb T}),$
 is the space of all $2\pi$-periodic generalized functions or the
  space of all $2\pi$-periodic Schwartz distributions. The space
 ${\mathcal E} ({\mathbb T}),$  can be characterized as follows: $F\in {\mathcal E} ({\mathbb T}),$
  if and only if it can be written as
$$
F( t )= \sum_{n=-\infty}^{\infty} F_n e^{i n t }\, , \quad
\mbox{where } \,F_n=\frac{1}{2\pi}<F,e^{-in t}>,
$$
and $$F_n=O\left( {n^p}\right), \quad as\quad |n| \rightarrow \infty,
\mbox{ for some } 0\leq p ,$$ or equivalently, $|F_n|\leq
B(|n|+1)^p$ for some $0<B, p$, and all $n$.  It can be shown that the action
of $F\in {\mathcal E} ({\mathbb T})$ on $f\in C^\infty ({\mathbb T})$ is given by
$$<F,f>=2\pi \sum_{n=-\infty}^\infty \overline{F}_n {\hat{f}}_n <\infty .$$

Let ${\mathcal A}$ be the space of all analytic functions on ${\mathbb T},$ or
equivalently
$${\mathcal A}=\left\{ f(t) =\sum_{n=-\infty}^{\infty} \hat{f}_n e^{i n t }, \quad \mbox{ with}\quad  \limsup_{|n|\rightarrow
\infty}\left|\hat{f}_n\right|^{1/|n|}<1 \right\}.$$ The topology of ${\mathcal A}$ is the topology it inherits
from $ C^\infty ({\mathbb T}) .$ Any $f\in
{\mathcal A}$ can be viewed as  the restriction of the harmonic
function
 $f(r,t) =\sum_{n=-\infty}^{\infty} \hat{f}_n r^{|n|}e^{i n t }$ to the unit circle.
 The dual space of ${\mathcal A},$ denoted by ${\mathcal A}^*,$
  is the space of hyperfunctions on  ${\mathbb T}$
 and can be characterized as
 $${\mathcal A}^*=\left\{ F(t) =\sum_{n=-\infty}^{\infty} F_n e^{i n t }, \quad \mbox{ with} \quad
\limsup_{|n|\rightarrow \infty}\left|F_n\right|^{1/|n|}\leq 1
\right\}.$$ A hyperfunction $F\in {\mathcal A}^*$ can be viewed as the
limit of the harmonic function $F(r,t)=\sum_{n=-\infty}^{\infty}
F_n r^{|n|}e^{i n t }$ as $r\rightarrow 1^-.$ The limit does not exist in the classical sense, but it exists in the sense
 of hyperfunctions.  Again the action of $F\in {\mathcal A}^*$ on $f\in {\mathcal A}$ is given by
$$<F,f>=2\pi \sum_{n=-\infty}^\infty \hat{f}_n \overline{F}_n <\infty .$$
The last series converges because $\limsup_{|n|\rightarrow
\infty}\left|F_n\hat{f}_n\right|^{1/|n|}<1. $ For more on hyperfunctions, see \cite{Johnson, Sato}.

Let ${\mathcal P}_{2\pi}$ be the space of all trigonometric polynomials
$\sum_{k=-n}^na_ke^{ikt}$ for all $n=0, 1, 2, \cdots ,$
endowed with the topology generated by the uniform convergence.
It has been shown that the dual space
${\mathcal P}^*_{2\pi}$ of ${\mathcal P}_{2\pi}$ is the space of all trigonometric
series $\sum_{k=-\infty}^\infty b_k e^{ikt}$
with no restriction on the growth of $b_k.$ This space was called the space of periodic
 generalized functions in \cite{Gorbacuk, Roumieu} and the space of  periodic ultra-distributions by Walter \cite{Walter}.

Now fix $0<q<1,$  and let $k\geq 1$ be a positive integer. Let
$$A_{q,k}=\left\{f(t): f(t)=\sum_{n=-\infty}^{\infty} \hat{f}_n
e^{i n t }, \quad \mbox{with } |\hat{f}_n|\leq c q^{|n|^k} \mbox{
for some } 0< c, \mbox{  and all } n\right\}.$$ It is easy to see that if $f\in
A_{q,k},$ then $|\hat{f}_n|\leq ce^{- \alpha |n|^k} ,$ for some
$\alpha , c>0, $ where $\alpha=\ln (1/q). $ Moreover,
 ${\mathcal A}=\cup_{0<q<1} A_{q,1},$ since, if $f\in A_{q,1},$ for some $0<q<1,$ then
 $\left|\hat{f}_n\right|\leq cq^{|n|}$ for some $0<c,$ and all $n;$ hence,
 $\left|\hat{f}_n\right|^{1/|n|} \leq c^{1/|n|} q,$ which implies that
 $\limsup \left|\hat{f}_n\right|^{1/|n|} <q<1.$ Therefore, $\cup_{0<q<1} A_{q,1} \subset {\mathcal A}.$
 Conversely, let $f\in {\mathcal A},$ then, $\limsup \left|\hat{f}_n\right|^{1/|n|}=q_1 <1$ for some $0<q_1<1.$ Therefore, for any $\epsilon>0$
 there exists $M$ such that $\left|\hat{f}_n\right|^{1/|n|}<(q_1+\epsilon)$ for all $M\leq |n|.$ Choose $\epsilon$ so that $(q_1+\epsilon)=q<1.$
 Let
 $$ C=\max_{|j|<M} \left\{1, \left|\hat{f}_j\right|/q^{|j|}\right\},$$
 then, it follows that
  $\left|\hat{f}_n\right|<Cq^{|n|}$ for all $n.$ That is $f\in A_{q,1}.$\\

Let $1<p$ such that $0<pq<1, $ and consider the space
$$A^{p,k}=\left\{F(t): F(t)=\sum_{n=-\infty}^{\infty} F_n e^{i n t }, \quad \mbox{with }
 |F_n|\leq c p^{|n|^k} \right\}.$$ It is easy to see that the
following inclusions hold
$${\mathcal P}_{2\pi}\subset A_{q,k}\subset A_{q,1}\subset {\mathcal A}\subset C^\infty ({\mathbb T}); \quad k\geq 1.$$
We provide $A_{q,k}$ with the topology  of uniform convergence. It
is easy to see that $A^{p,k}$ is a subset of the dual space of
$A_{q,k},$ and we have the inclusions
$${\mathcal E} ({\mathbb T}) \subset {\mathcal A}^*\subset A^*_{q,1}\subset A^*_{q,k}\subset {\mathcal P}^*_{2\pi}.$$
Moreover,
for $F\in A^{p,k}$ and $f\in A_{q,k}$ we have $<F,f>=2\pi
\sum_{n=-\infty}^{\infty} \hat{f}_n \overline{F_n}<\infty $ and
the last series converges because
$$\left| F_n\hat{f}_n\right|\leq c (pq)^{|n|^k}\leq c(pq)^{|n|}$$ and the
series $\sum_{n=-\infty}^\infty (pq)^{|n|}$ converges because $0<pq<1.$
The dual of $A_{q,k}$ is the countable-union space $A^*_{q,k}
=\cup_{p<1/q} A^{p,k}.$

Recall that if $A$ is a space of functions defined on some interval $I,$ then $A$
is said to be quasi-analytic if for any $f\in A,$ with
$f^{(n)}(x_0)=0$ for all $n=0,1, 2, \cdots $ and $x_0\in I,$ then
$f$ is identically zero. A class of quasi-analytic functions
cannot contain any function that vanishes on an interval of
positive measure.

\begin{lemma}
The space $A_{q,k}$ is a space of quasi-analytic functions.
\end{lemma}
{\bf Proof:}
First, we show that if $f\in A_{q,k},$ then there exist $0< c_f
,B$ such that
$$\left| f^{(m)}(x)\right| \leq c_f B^m m^{m/k}.$$
Since
$$f^{(m)}(x)=\sum_{n=-\infty}^\infty (in)^m \hat{f}(n) e^{inx},$$
we have
\begin{eqnarray*}
\left| f^{(m)}(x)\right| & \leq & \tilde{c} \sum_{n=-\infty}^\infty |n|^m
q^{|n|^k} \leq 2c \int_0^\infty x^m e^{-ax^k} dx\\
&=& \frac{2c}{k}\left( \frac1a\right)^{(m+1)/k}\int_0^\infty
y^{(m+1)/k-1} e^{-y} dy=CB^m \Gamma\left(\frac{m+1}{k}\right)\\
&=& CB^m \Gamma\left(\frac{m+1-k}{k}+1\right),
\end{eqnarray*}
where
$$a=\ln \frac{1}{q}>0, \; B=\left(\frac{1}{a}\right)^{1/k}, \;
C=\frac{2c}{k}\left( \frac1a\right)^{1/k}.$$ But by Sterling's
formula $\Gamma (x+1)\approx \sqrt{2\pi x}e^{-x} x^x,$ we obtain after some calculations
$$\left| f^{(m)}(x)\right|\leq C_1B_1^m m^{m/k}, \quad \mbox{ for
some } C_1 , B_1 .$$

Following Rudin's notation \cite[p. 410-416]{Rudin}, we have for any $f\in A_{q,k}$
$$\left| f^{(m)}(x)\right| \leq c_f B^m M_m \mbox{ where } M_m= m^{m/k}.$$
but
$$\sum_{m=1}^\infty \frac{1}{M^{1/m}_m}=\sum_{m=1}^\infty
\frac{1}{m^{1/k}}=\infty, \quad \mbox{since} \quad k \geq 1 ,$$
it follows from the Denjoy-Carleman theorem that
the space $A_{q,k}$ is quasi-analytic. \ep

\begin{defi}
 A generalized function $F\in {\mathcal P}^*_{2\pi}$ is said to be positive
if $f\in {\mathcal P}_{2\pi}$ is positive, then $<F,f> $ is positive.
\end{defi}

Now we prove our main result in this section, which is the analog of Theorem \ref{thetaprop},
\begin{theorem}
For each $0\leq t$ and $F\in A^*_{q,k},$ define the multiplier
operator ${\mathcal T}_t$ by
$${\mathcal T}_t F(x)=\sum_{n=-\infty}^\infty F_ne^{-n^2 t} e^{inx}.$$
Then
\begin{enumerate}
\item[i)] ${\mathcal T}_t F(x)$ converges weakly in $A^*_{q,k},$ for all $t\geq
0$
 \item[ii)] ${\mathcal T}_{t_1+t_2}={\mathcal T}_{t_1}{\mathcal T}_{t_2}$ and $T_0=I. $
  \item[iii)]  ${\mathcal T}_tF \rightarrow F$ weakly as $t\rightarrow 0.$
  \item[iv)] The limit
$$\lim_{t\rightarrow 0}\frac{{\mathcal T}_t F-F}{t}=\frac{d^2 F}{dx^2},$$
converges weakly, where the derivative on the right-hand side is taken in
the generalized function sense, i.e.,
$$<F^{(n)},f>=(-1)^n <F, f^{(n)}>,$$
and $f\in A_{q,k}$.
 \item[v)] $u(x,t)  =  {\mathcal T}_t F(x)$ is a weak solution of the heat equation on ${\mathbb T} \times\mathbb{R}^+$,
  $$ \frac {\partial  u(x,t) }{\partial t}  =\frac {\partial^2 u(x,t) }{\partial x^2}, \quad (x,t) \in {\mathbb T} \times\mathbb{R}^+,$$
  with initial condition $u(x,0)  = F(x), \, x \in{\mathbb T},$
\item[vi)] For every $F\in A^*_{q,2}$
there exists $t_F$ so that ${\mathcal T}_tF \in A_{p,2}$ for all $t\geq t_F.$
\item[vii)] $T_t$ is a positive operator on $A^*_{q,k},$ meaning ${\mathcal T}_t F\geq 0$ if $F\geq
0.$

\end{enumerate}
\end{theorem}
 \dem
   i)  Since $F\in A^*_{q,k},$ we have $F\in A^{p,k}$ for some $1<p.$ Thus, for all $f \in A_{q,k}$
 $$\langle {\mathcal T}_tF, f\rangle = 2\pi \sum_{n=-\infty}^\infty
 \hat{f}_n  \overline{F_n}e^{-n^2t} $$
 and we have  for all  $ 0\leq t ,$
\begin{eqnarray*}
\left| \langle {\mathcal T}_tF, f\rangle \right| & \leq & 2\pi
\sum_{n=-\infty}^\infty \hat{f}_n \overline{F_n}e^{-n^2t} \leq
2\pi \sum_{n=-\infty}^\infty
 \left|F_n\right| \left|\hat{f}_n\right|\\
 &\leq & \tilde{C} \sum_{n=-\infty}^\infty (pq)^{|n|^k} <C \sum_{n=-\infty}^\infty (pq)^{|n|} <
 \infty.
 \end{eqnarray*}
 The last series converges because $0<(pq)<1. $

 ii) This is trivial.

 iii) First
 $$\left| \langle {\mathcal T}_tF-F, f\rangle \right| \leq 2\pi \sum_{n=-\infty}^\infty
 \left|F_n\right| \left|\hat{f}_n\right|\left( e^{-n^2t}-1\right) <\infty .$$
 Then, as in part i) since
 $$ \sum_{n=-\infty}^\infty
 \left|F_n\right| \left|\hat{f}_n\right| <\infty,$$
 we can take the limit as $t\rightarrow 0$ inside the summation to
 obtain  ${\mathcal T}_tF
\rightarrow F$ weakly as $t\rightarrow 0$.

 iv) Let $f\in A_{q,k}.$ Then we have
  $$\langle \frac{{\mathcal T}_tF-F}{t}, f\rangle = 2\pi\sum_{n=-\infty}^\infty
 \hat{f}_n \overline{F_n } \left( \frac{e^{-n^2t}-1}{t}\right)  .$$
By taking the limit as $t\rightarrow 0,$ we obtain
$$\lim_{t\rightarrow 0} \langle \frac{{\mathcal T}_tF-F}{t}, f\rangle = 2\pi\sum_{n=-\infty}^\infty
 (-n^2)\hat{f}_n \overline{F_n } .$$
 But on the other hand,
 $$ \langle \frac{d^2 F}{dx^2}, f\rangle  =\langle F, \frac{d^2 f}{dx^2}\rangle
= 2\pi\sum_{n=-\infty}^\infty
 (-n^2)\hat{f}_n \overline{F_n }.$$
 Thus,
$$\langle \frac{{\mathcal T}_tF-F}{t}, f\rangle =\langle \frac{d^2 F}{dx^2}, f\rangle
,$$ which implies iv).

v) The proof follows from (i) and (iv).

 vi) Let $F\in A^*_{q,2},$ then there exists $1\leq p$ such that
$F\in A^{p,2}$. Choose $t_F$ so that $p<e^{ t_F /2}$. Hence, for
$t_F\leq t ,$ we have
$${\mathcal T}_tF(x)=\sum_{n=-\infty}^\infty F_n e^{-n^2 t}e^{inx}$$
and
\begin{eqnarray*}
\left| {\mathcal T}_tF(x)\right|  & \leq & \sum_{n=-\infty}^\infty
\left|F_n\right| e^{- n^2 t}
 \leq   c \sum_{n=-\infty}^\infty p^{n^2}e^{-n^2 t_F}\\
&=& c \sum_{n=-\infty}^\infty \left(
\frac{p}{e^{t_F}}\right)^{n^2}<\infty
\end{eqnarray*}

Moreover, $$\left| F_n\right|e^{-n^2t}\leq c p^{n^2} e^{-n^2t}\leq
c p^{n^2} e^{-n^2t_F}\leq c e^{n^2 t_F /2} e^{-n^2t_F}\leq c
e^{-n^2 t_F /2}\leq c q^{n^2},$$ where $q=e^{-t_F/2}<1 .$ Thus,
$T_tF(x)\in A_{q,2}.$

vii) Let $f\in A_{q,k},$ then $(f\star\theta_3)(x)\in A_{q,k}$. For
\begin{eqnarray*}
(f\star\theta_3)(x)&=&
\frac{1}{2\pi}\int_0^{2\pi}\left(\sum_{m=-\infty}^\infty\hat{f}_m
e^{imy}\right) \left(\sum_{n=-\infty}^\infty e^{-n^2t}e^{in(x-y)}\right)dy\\
&=&\frac{1}{2\pi}\sum_{m,n=-\infty}^\infty\hat{f}_m
e^{-n^2t}\int_0^{2\pi}e^{inx}e^{iy(m-n)}dy\\
&=& \sum_{n=-\infty}^\infty e^{-n^2t}\hat{f}_n e^{inx}.
\end{eqnarray*}
But since $\left| \hat{f}_n e^{-n^2t}\right|\leq \left| \hat{f}_n
\right|\leq cq^{n^2},$ it follows that $(f\star\theta_3) \in
A_{q,k}.$

Now let  $F\in A^*_{q,k}$
\begin{eqnarray}
\langle {\mathcal T}_tF,f\rangle &=& 2\pi \sum_{n=-\infty}^\infty e^{-n^2t}
\overline{F}_n \hat{f}_n=  \sum_{n=-\infty}^\infty e^{-n^2t}
\langle
F, e^{iny}\rangle {\hat{f}}_n \nonumber\\
&=& \frac{1}{2\pi} \langle F,  \sum_{n=-\infty}^\infty e^{-n^2t}
e^{iny}\langle
f(z) ,e^{-inz} \rangle \rangle \nonumber \\
&=& \frac{1}{2\pi} \langle F,  \langle f, \sum_{n=-\infty}^\infty
e^{-n^2t} e^{in(y-z)} \rangle \rangle . \label{star1}
\end{eqnarray}
But
$$\langle f, \sum_{n=-\infty}^\infty e^{-n^2t}
e^{in(y-z)} \rangle = \langle f , \theta_3(y-., e^{-t}) \rangle
=2\pi (f\star\theta_3 ).$$ Therefore, from (\ref{star1}),  we have
$$\langle {\mathcal T}_tF,f\rangle =\langle F,f\star\theta_3(y-., e^{-t})
\rangle ,$$ which is well defined because $(f\star\theta_3) \in
A_{q,k},$ and $F\in A^*_{q,k}.$

 Finally, we need to show that ${\mathcal T}_t F$ is positive if $F$ is positive. Let $f$
be positive and since $\theta_3$ is positive by Proposition 2.1,
it follows that $f\star\theta_3$ is positive. Since $F$ is
positive, it follows that $\langle F,f\star\theta_3(y-., e^{-t})
\rangle$ is positive which in turn implies that $\langle
{\mathcal T}_tF,f\rangle $ is positive for any positive function $f\in
A_{q,k},$ and this completes the proof. \ep

\end{document}